\documentclass[11pt]{article}
\usepackage{amsmath}
\usepackage{amsfonts}
\usepackage{amssymb}
\usepackage{amsthm}

\newcommand{\Z}{{\bf{Z}}}
\newcommand{\Q}{{\bf{Q}}}

\newcommand{\R}{{\bf{R}}}
\newcommand{\C}{{\bf{C}}}

\newcommand{\T}{{\bf{T}}}
\newcommand{\G}{{\bf{G}}}

\newcommand{\A}{{\bf{A}}}

\newcommand{\pp}{{\mathfrak{p}}}
\newcommand{\tensor}{\otimes}
\newcommand{\cross}{\times}
\newcommand{\ra}{\rightarrow}

\newcommand{\OO}{\mathcal{O}}

\newcommand{\Hc}{{{\mathcal{H}}}}
\newcommand{\Ht}{{{\mathfrak{H}}_{\scriptscriptstyle{2}}}}

\newcommand{\rram}{{r_{\scriptscriptstyle{B}}}}

\newcommand{\rinert}{{r_{\scriptscriptstyle{K}}}}

\newcommand{\fram}{{f_{\scriptscriptstyle{B}}}}
\newcommand{\finert}{{f_{\scriptscriptstyle{K}}}}

\newcommand{\p}{{\mathfrak{p}}}
\newcommand{\n}{{\mathfrak{n}}}
\newcommand{\af}{{\mathfrak{a}}}

\newcommand{\Jtau}{{J_\tau}}
\newcommand{\Hom}{{\rm Hom}}
\newcommand{\rec}{{\rm rec}}

\newcommand{\Pic}{{\rm Pic}}

\newtheorem{lem}{Lemma}[section]

\newtheorem{conj}[lem]{Conjecture}
\newtheorem{thm}[lem]{Theorem}

\theoremstyle{definition}

\newtheorem{rmk}[lem]{Remark}

\newcommand{\later}[1]{}
\newcommand{\comment}[1]{}
\newcommand{\com}[1]{}

\newcommand{\noi}{\noindent}

\newcommand{\thetitle}
{Darmon points on elliptic curves over totally real fields}

\begin{document}
\parindent=2em

\date{}

\title{\thetitle}
\author{Amod Agashe and Mak Trifkovi\'c}
\maketitle


\com{need orientation only to get free action of PicO. ?? }

\section{Introduction}
Let $F$ be a number field and let $E$ be an elliptic curve over~$F$ of
conductor an ideal~$N$ of~$F$.  We assume throughout that $F$ is
totally real: in that case, it is known, under minor hypotheses, that
there is a newform~$f$ of weight~$2$ on~$\Gamma_0(N)$
over~$F$\later{cite weil} whose $L$-function coincides with that
of~$E$ 
(see, e.g., \cite{zhang-annals} or \cite[\S7.4]{dar-rat}).
We fix a quadratic
extension $K/F$.  When $K$ is totally complex, a classical theory
produces a family of Heegner points of $E$, defined over ring class
fields of $K$.  The Galois action on them is given by a Shimura
reciprocity law, and their heights relate to the derivative
of the $L$-function of~$E$ over~$K$ at~$1$
(see, e.g., \cite{zhang-annals} or \cite[\S7.5]{dar-rat}).

We assume therefore from now on that $K$ has at least one real place.
The goal of the theory of Darmon points
(earlier called Stark-Heegner points) is to extend to such $K$ the  
construction of Heegner points.  In ~\cite{dar-annals}, Darmon started the theory in the case
$F=\Q$ and $K$ a totally real quadratic field.  
In~\cite{greenberg}, Greenberg generalizes this work
to give a conjectural construction of points
in the case where $F$ is arbitrary totally real quadratic field of narrow class number one, $E$ is semistable, $N \neq (1)$,
the sign of the functional equation of~$E$ over~$K$ is~$-1$,
there is a prime dividing the conductor of~$E$ that is inert in~$K$,
and the discriminant of~$K$ is coprime to~$N$. The techniques
used are $p$-adic in nature.
In~\cite{gartner}, Gartner generalizes the work of~\cite[Chap~8]{dar-rat} and~\cite{darlog}
to give a construction of what he calls Darmon points in certain situations using
archimedian techniques. The conditions under which Gartner's construction works
are a bit too technical to describe here, but we shall describe them 
in Section~\ref{sec:gartner}.   All the constructions share a basic
outline: one computes an archimedean or $p$-adic integral of the
modular form associated to $E/F$, and plugs the resulting value into a
Weierstrass of Tate parametrization of $E$ to produce the Darmon point.


Let $\OO \subseteq K$ be an~$\OO_F$-order such that 
${\rm Disc}(\OO/\OO_F)$ is coprime to~$N$. 
In this article, 
we show that if the sign of the functional equation of~$E$ over~$K$ is~$-1$,
the discriminant of~$K$ is coprime to~$N$,
and the part of~$N$
divisible by primes that are inert in~$K$ is square-free, then
one can apply either the construction of Gartner or the construction of Greenberg 
(after removing the assumption that $F$ has narrow class number one, which we show how to do) to conjecturally associate to~$\OO$ a point 
that we call a Darmon point (actually there
are choices in the construction, and so one gets several Darmon points).
This point is intially defined over a transcendental extension of~$K$, but
we conjecture that the point is algebraic, defined over the narrow ring
class field extension of~$K$ associated to the order~$\OO$. 
This point comes with an action of the narrow class group of~$\OO$,
and we state a conjectural Shimura receprocity law for this action.

In Section~\ref{sec:gg} 
we recall and slightly modify the constructions of Greenberg and Gartner.
In Section~\ref{shdetails} we show how one of the two constructions
can be carried out under our hypotheses.
We assume throughout this article
that the reader is familiar with~\cite{greenberg} and~\cite{gartner}.

\section{The constructions of Gartner and Greenberg} \label{sec:gg}

In this section we discuss Gartner's and Greenberg's constructions. 
Gartner makes several assumptions that are sometimes not
made explicit in~\cite{gartner};
we clarify what hypotheses are needed in Gartner's construction
and also modify it a bit
so that it can be unified better with Greenberg's construction.
We also show how to generalize
Greenberg's construction to remove the class number one assumption in~\cite{greenberg}.

Both constructions require the existence of a suitable
quaternion algebra~$B$ 
in order to use the Jacquet-Langlands
correspondence. So let $B$ be a
quaternion algebra over~$F$. 
We will impose certain assumptions that $B$ (and other objects) will
have to satisfy in each of Gartner's or Greenberg's constructions.
In Section~\ref{shdetails} we shall explain when these assumptions are met.

First, in either construction, one  needs an embedding of~$K$ into~$B$.
Recall that $B$ is said to be split at a place $v$ of~$F$ if
$B \tensor_F F_v$ is the matrix algebra, and ramified at $v$ otherwise.
It is known that a quaternion algebra is determined up to
isomorphism by the set of ramified places, which is
finite of even cardinality.   Conversely, for any finite set of
places of
even cardinality there is a quaternion algebra ramifying at
these places.
We say that a real place of~$F$ splits in~$K$ if there are two real places
of~$K$ lying over it,
and we say that it is inert 
\later{go back to ramified since do so for finite primes?}
otherwise (such a place 
is usually said to be ramified, but we prefer
to call it inert 
to avoid confusing ramification in~$K$ with ramification
of~$B$).
\later{check: even for finite primes?}

\ \\
\noi {\it Assumption A:} Assume that there is an embedding
of~$q:K\hookrightarrow B$,
i.e., that each place where $B$ ramifies, archimedean or not, is inert in~$K$. 
\\


\subsection{Gartner's construction} \label{sec:gartner}

We now outline the construction of Gartner, along with some modifications.
For details and proofs of the claims made below, please
see~\cite{gartner}. We try to use notation consistent with
or similar to that in~\cite{gartner} as much as possible.

We start by listing the assumptions used in Gartner's construction.
Let $d$ denote the degree of~$F$ over~$\Q$ and let $\tau_1, \ldots, \tau_d$
denote the archimedian places of~$F$. \\

\noi {\it Assumption B1: }Suppose that there is exactly one 
archimedian place of~$F$ where $B$ is split but which does
not split in~$K$.  \\
\later{gives Weier unif. Split into 2? (exists and unique?)
and second unnecessary in Gartner? what if $r_K \neq 0, 1$?}

Without loss of generality, assume that the
archimedian place of~$F$ where $B$ is split but which does
not split in~$K$ is~$\tau_1$. 
Let $r$ be the integer such that the archimedian places of~$F$
that split in~$K$ are $\tau_2, \ldots, \tau_r$; since
$K$ is not a CM field,
$r \geq 2$. By our Assumption~A,
$B$ necessarily splits at $\tau_1, \ldots, \tau_r$ and
by Assumption~B1, it necessarily
ramifies at~$\tau_{r+1}, \ldots, \tau_d$.

If $S$ is a ring, then let $\widehat{S}$ denote $S\tensor_\Z \widehat{\Z}$.
Let $R$ be an Eichler order of~$B$.\\

\noi {\it Assumption B2:} Assume that $f$ corresponds to an automorphic form 
on~$\widehat{R}$ \later{check}
under the Jacquet-Langlands correspondence.  \\

Let
$b \in \widehat{B}^\times$.
Let $\OO \subseteq K$ be an~$\OO_F$-order such that 
${\rm Disc}(\OO/\OO_F)$ is coprime to~$N$. In order to get a Darmon
point in the narrow class field 
associated to the order~$\OO$, along  
with an action of~$\Pic(\OO)^+$, we make the following
assumption (which is not made in~\cite{gartner}):\\
\later{needed for action or just for free action? check Mak's article}

\noi {\it Assumption B3:} 
Suppose that $q(K) \cap b \widehat{R}b^{-1} = q(\OO)$,
i.e., that $q$ is an optimal embedding of~$\OO$ into
the order~$B \cap b \widehat{R}b^{-1}$. \\

We now start the construction. 
Let $G = {\rm Res}_{F/\Q} B^\cross$ and 
let $\A_f$ denote the set of finite adeles over~$\Q$.
Let $H = \widehat{R}^\times$ and
let ${\rm Sh}_H(G)$  
denote the quaternionic Shimura variety \com{over a certain reflex field}
whose complex valued points are given by
$G(\Q) \backslash (\C \setminus \R)^r \times G(\A_f)/H$.
Let $b \in \widehat{B}^\times$. 
Let $T = {\rm Res}_{K/\Q}(\G_m)$. 
The embedding~$q$ induces an embedding
of~$T$ in~$G$ that we again denote by~$q$ for simplicity. 
Using the embeddings associated to~$\tau_1, \ldots, \tau_r$,
where $B$ splits, we get a natural action of 
$q(T(\R)^0)$ on~$X = (\C \setminus \R)^r$.
Let $T^0$ be a fixed orbit
of~$q(T(\R)^0)$ whose projection to the first component of~$X$ 
is a point (recall that~$\tau_1$ is a complex place); we fix this point henceforth and 
denote it by~$z_1$.
Let $T_{q,b}$
denote the projection to~${\rm Sh_H}(G)(\C) = G(\Q) \backslash (\C \setminus \R)^r \times G(\A_f)/H$
of~$T^0
\times G(\A_f)$.
Then $T_{q,b}$ is a torus of
dimension $r-1$.
Note that our torus corresponds to the torus denoted
${\mathcal T}^0_b$ 
in Section~4.2 of~\cite{gartner}; moreover Gartner actually works
with a modified Shimura variety denoted ${\rm Sh}_H(G/Z, X)$ in loc. cit. 
However, the construction
goes through mutatis mutandis with ${\rm Sh}_H(G)$ as well, which is what we shall do in this article.
Thus our construction is a slightly modified version of that of Gartner. 
We are doing this modification
to get a version of the (conjectural) Shimura receprocity law that
is similar to that of Greenberg (Conjecture~3 in~\cite{greenberg}).
Using the theorem of Matsushima and Shimura~\cite{matshi}, one shows that there
is an $r$-chain that we denote~$\Delta_{q,b}$ (called~$\Delta_b$ in~\cite{gartner})
on~${\rm Sh}_H(G)$ whose boundary is an integral multiple of~$T_{q,b}$. 
This uses Proposition~4.5 in~\cite{gartner}, whose proof assumes that the Shimura variety is compact, i.e., that $B$ is not
the matrix algebra. If $B$ is the matrix algebra, then 
we are in the ATR setting dealt with in~\cite[Chap.~8]{dar-rat} 
and~\cite{darlog}; 
thus in this article, we are subsuming the ATR construction under
Gartner's construction
(in fact, Gartner's work was motivated by the ATR construction).

Let $\phi$
denote the automorphic form on~$H$ corresponding to~$f$ under
the Jacquet-Langlands correspondence (recall our Assumption~B2).
Analogous to the construction of the form denoted~$\omega_\phi^\beta$
in~\cite{gartner}, we get
a form that we denote~$\omega_\phi$ 
on~${\rm Sh}_H(G)$ by taking $\beta$ to be the trivial character
(we could allow $\beta$ to be an aribitrary character, but we are taking it to be
the trivial character for the sake of simplicity and also to get an action
of the narrow class group below). \com{Else we get an action of the  subgroup
$K_A^\cross / ( K_+^\cross (K \tensor R)_+^\cross \widehat{O}^\cross$. )}
Assuming conjectures
of Yoshida \cite{yoshida}, the periods of~$\omega_\phi$ form a lattice that is homothetic to
a sublattice of the Neron lattice~$\Lambda_E$ of~$E$. 
Then the image of a suitable integer multiple of~$\int_{\Delta_{q,b}} \omega_\phi$ 
is independent of
the choice of the chain~$\Delta_{q,b}$ made above. 
Let $\Phi: \C/\Lambda_E \ra E(\C)$ denote
the Weierstrass uniformization of~$E$. Then the Darmon point~$P_{q,b}$ in~$E(\C)$ is
defined as a suitable multiple of the image of~$\int_{\Delta_{q,b}} \omega_\phi$ 
under~$\Phi$ (our point corresponds to the point~$P^\beta_b$
in~\cite{gartner}). 
It is conjectured that the point $P_{q,b}$ in~$E(\C)$ has algebraic coordinates.

Let $\widehat{q}:\hat K\to\hat B$ denote the map obtained from~$q$ by tensoring with~$\widehat{F}$.
Let $K_\A$ denote the ring of ad\`{e}les of~$K$.  Denote by $a_f$ the
non-archimedean part of $a \in K_\A$.
Following Gartner, we define an action of $K_\A^\times$ on Darmon points~$P_{q,b}$ by
$$a * P_{q,b} = P_{(q, \widehat{q}(a_f) b)}.$$ 
An easy check shows that the new pair satisfies Assumption~B3.

Denote by $K_+$ the subset of elements of~$K$ that are positive in all
real embeddings.
As usual, $K^\times$ is
embedded into $K_\A^\times$ diagonally.
We claim that the action above factors through
$\widehat{\OO}^\times (K \tensor_\Q \R)^\times K_+^{\times}$, 
i.e., we have an action of
$\widehat{\OO}^\times \backslash K_\A^\times /(K \tensor_\Q \R)^\times K_+^{\times}$, which
is the narrow class group~$\Pic(\OO)^+$. 
To prove the invariance under the action of~$\widehat{\OO}$, 
note that by our condition above, $\widehat{q}(\widehat{\OO}) \subseteq b \widehat{R}^\times b^{-1}$,
so if $a \in \widehat{\OO}$, then $\widehat{q}(a_f) = b r b^{-1}$ for some $r \in \widehat{R}$. Hence
$\widehat{q}(a_f) b = b r b^{-1} b = b r$, and so 
 $a * P_{q,b} = P_{(q, br )} = P_{q,b}$ since $\widehat{R}$ acts trivially
on ${\rm Sh}_H(G)$. 
Next $(K \tensor_\Q \R)^\times$ clearly acts trivially.
It remains to show invariance under the action of~$K_+^\times$.
If $k \in K_+^\times$, then 
 $k * P_{q,b} = P_{(q, \widehat{q}(k) b)}$. 
Let $x_0 \in T^0$.
Then $T_{(q, q(k) b)}$ consists of
images of points of the form
$(y, q(k) b)$ in~$(\C \setminus \R)^r \times G(\A_f)$
such that $y = t x_0$ for some $t \in q(T(\R)^0)$. 
Letting $\pi$ denote the projection map from $(\C \setminus \R)^r \times G(\A_f)$ 
to~${\rm Sh}_H(G)(\C) = G(\Q) \backslash (\C \setminus \R)^r \times G(\A_f)/H$, we have
$\pi(y, q(k) b) = \pi((t x_0, q(k) b) = \pi(q(k^{-1}) t x_0, b)
= \pi(t q(k^{-1}) x_0, b) $, as elements of~$q(K)$ commute. 
Thus the point $P_{(q, \widehat{q}(k) b)}$ is obtained from the orbit with
base point~$q(k^{-1}) x_0$. 
Recall that the projection of~$q(T(\R)^0)$ to the first component of~$X$ 
is a point~$z_1$, and thus $z_1$ is
fixed by $q(K)$. In particular, the projection of the orbit of the point $q(k^{-1}) x_0$ 
to the first component of~$X$ is again~$z_1$. Moreover the projections of
the orbits of the point $q(k^{-1}) x_0$ and of the point~$x_0$ lie in the same
connected component of each copy of~$\C \setminus \R$ in~$X$ since $k \in K_+$.
In view of the last two sentences in~\cite[Prop~4.7]{gartner}, 
$P_{(q, \widehat{q}(k) b)} = P_{(q, b)}$, which finishes our proof
of the claim.

Thus we get an action of
the narrow class group~$\Pic(\OO)^+$ 
on Darmon points~$P_{q,b}$; we denote this action again by~$*$. 
\com{Comment: this action may be given more
explicitly in terms of ideals as in Gross.}
\begin{conj} \label{conj:gartner}
The point~$P_{q,b}$ is defined over the narrow ring
class field extension~$K_\OO^+$ of~$K$ associated to the order~$\OO$.
If $\alpha \in \Pic(\OO)^+$, then 
$\rec(\alpha)(P_{q,b}) = \alpha * P_{q,b}$, where 
$\rec: \Pic^+ \OO \ra {\rm Gal}(K_\OO^+/K)$ is the reciprocity
isomorphism of class field theory. 
\end{conj}
\com{Comment: This implies (as in Gartner 3.1.3.1) that the point $P_{q,b}$ is defined 
over the narrow ring class field of~$\OO$ (under our constraints on~$(q,b)$ -- see p. 55 and 56
of Gartner's thesis;
$q_{\A}^{-1}(b H b^{-1} B_+^\times) = \widehat{\OO} (K \tensor_\Q \R)^\times_+ K^{\times}$?), 
which is denoted $K^+_b$ in Gartner.}

The following assumption is not needed for Gartner's construction,
but we shall mention it since it will be useful in Section~\ref{shdetails}
(see also Remark~\ref{rmk:assn}(i)).\\

\noi {\it Assumption B4:} 
Suppose that the finite primes where
$B$ is ramified are exactly the 
primes that divide~$N$ and are inert in~$K$. \\


\later{Gartner assumes B \neq M_2, needed so Mat Shi suffices. Use ATR}

\subsection{Greenberg's construction} \label{sec:greenberg}  

We now discuss the construction of Greenberg, and
show how to remove assumption that~$F$ has class number one made in~\cite{greenberg}.
For details of the construction, please
see~\cite{greenberg}.   We start by listing the assumptions needed.
Recall that $d$ denotes the degree of~$F$ over~$\Q$ and $\tau_1, \ldots, \tau_d$
denote the archimedian places of~$F$. 
Since $K$ is not CM, there is at least one 
infinite place of~$F$ that splits in~$K$.
Let $n$ denote the number of such places,
and without loss of generality, assume that these places $\tau_1,
\ldots, \tau_n$ (in the previous section, we wrote $r$ instead of $n$;
we change notation to be consistent with
or similar to that in~\cite{greenberg} as much as possible).
By Assumption~A, $B$ is split at $\tau_1, \ldots, \tau_n$ and
can be ramified or split at $\tau_{n+1}, \ldots, \tau_d$.
However we insist:\\

\noi{\it Assumption C1:} Suppose that $\tau_{n+1}, \ldots, \tau_d$ are
precisely the infinite primes where
$B$ ramifies.\\
\later{rK = 0}

\noi{\it Assumption C2:}
Suppose that there is a prime ideal~$\p$ of~$F$ that exactly divides~$N$ 
and is inert in~$K$. \\
\later{fK = 1; needed for tate param; fK > 1 OK?}


\noi{\it Assumption C3:} Suppose that the part of~$N$
divisible by primes that are inert in~$K$ is square-free and that the 
finite primes where
$B$ is ramified are exactly the 
primes other than~$\p$ that divide~$N$ and are inert in~$K$. \\
\later{replace by JL}

Let $\n$ be the part of~$N$ supported at primes other than~$\p$ that divide~$N$ 
and where $B$ is split. For each ideal~$\af$ of~$\OO_F$ coprime to the
discriminant of~$B$, choose an Eichler order~$R_0(\af)$ in~$B$ of level~$\af$
as in~\cite[\S2]{greenberg}.
Let $R = R_0(\n)$ be the Eichler order in $B$ of level $\n$. 
As in Section~\ref{sec:gartner}, 
we choose 
$b \in \widehat{B}^\times$ and 
impose the analog of Assumption~B3:\\

\noi {\it Assumption C4:} 
Suppose that $q(K) \cap b \widehat{R} b^{-1} = q(\OO)$,
i.e., that $q$ is an optimal embedding of ~$\OO$ into
the order~$B \cap b \widehat{R}  b^{-1}$.\\
\later{needed to get action and Sh rec or to get phs? Sim for B3}

We remark that the assumptions made above are not exactly the assumptions
made in~\cite{greenberg}, but suffice for the construction (e.g., the assumption
made in~\cite{greenberg} that the sign of the functional equation of~$E$
over~$K$ is~$-1$ is used to show that a quaternion algebra~$B$ satisfying Assumptions~C1
and~C3 exists).

We now describe the construction. As in Section~\ref{sec:gartner}, let
$G = {\rm Res}_{F/\Q} B^\cross$, let $H = \widehat{R}^\times$,
and let ${\rm Sh}_H(G)$ 
denote the quaternionic Shimura variety \com{over a certain reflex field}
whose complex points are given by
$G(\Q) \backslash (\C \setminus \R)^n \times G(\A_f)/H$, where
$n$ 
is the number of real places of~$F$ where $B$ splits
($n$ is the same as the~$r$ in Section~\ref{sec:gartner}).
Let $G(\R)^0$ denote the identity component of~$G(\R)$ and let $G(\Q)^0 = G(\R)^0 \cap G(\Q)$.
Let $C \subseteq \widehat{B}^\times$ be a system of representatives of the double cosets
$G(\Q)^0 \backslash \widehat{B}^\times / H$. If $g \in \widehat{B}^\times$, then let $\Gamma'_g
= gHg^{-1} \cap G(\Q)^0 \subseteq G(\R)^0$ and let $\Gamma_g$ denote
the natural projection of the image of~$\Gamma'_g$ in~$PGL_2^+(\R)^n$ 
(the projection is obtained
via the embeddings associated to the places~$\tau_1, \ldots, \tau_n$
where $B$ splits). Let $\Ht$ denote the upper half plane.
Then ${\rm Sh}_H(G)(\C)$ 
is homeomorphic to the disjoint union of $\Gamma_g \backslash (\Ht)^n$
as $g$ ranges over elements of~$C$. 
%
%
Greenberg assumes that the narrow class number of $F$ (and therefore
of $B$) is one, in which case $C$ is a singleton set and 
$\Gamma_g = \Gamma_0(\n)$.  
Greenberg's construction uses group homology and cohomology for the 
group~$\Gamma_0(\n)$ with coefficients in various modules. 
When the narrow class number is not one, one has to replace 
the homology groups of~$\Gamma_0(\n)$ with
the direct sum over $g \in C$ of the homology groups of~$\Gamma_g$.

As in Section~\ref{sec:gartner}, we 
construct the torus $T_{q,b}$ in~${\rm Sh}_H(G)$. 
The inclusion of the torus in~${\rm Sh}_H(G)$ induces a map on the corresponding
$n$-th homology groups.
The image under this map of a generator of the $n$-th homology group of the torus
gives an element of the 
$n$-th homology group of~${\rm Sh}_H(G)$
(in fact, since the torus is connected, that element lies in one of
direct summands in the homomlogy); this
element replaces the element denoted~$\Delta_\psi$ starting with Lemma 21
in~\cite{greenberg} (note that in loc. cit.,
before Lemma 21 , $\Delta_\psi$ is considered
to be an element of a homology group of~$\Gamma_0(\n)_\psi$, but starting with
Lemma 21, $\Delta_\psi$ is considered
to be an element of a homology group of~$\Gamma_0(\n)$).
Greenberg also uses homology groups of~$\Gamma_0(\n)$ with coefficients
in a module denoted~${\rm Div\ } \Hc_\p^\OO$ in loc. cit.; here,
$\Hc_\p^\OO$ denotes a certain set of points, one corresponding
to each optimal embedding of~$\OO$ in~$R$
(see page~561 of loc. cit. for details).
To generalize this construction, for each $g \in C$, we
define $\Hc_{\p,g}^\OO$ to be the analogous set of points, which is in bijection
with the set of optimal embeddinga of~$\OO$ in~$B \cap g \widehat{R}g^{-1}$.
Then the homology groups $H_i(\Gamma_0(\n), {\rm Div\ }\Hc_\p^\OO)$
get replaced by $\oplus_{g \in C} H_i(\Gamma_g, {\rm Div\ }\Hc_{\p,g}^\OO)$.
On the group cohomology side, Greenberg considers cohomology groups of 
the groups~$\Gamma_0(\n)$ and~$\Gamma_0(\p \n)$. The cohomology groups
of~$\Gamma_0(\n)$ again get replaced by 
the direct sum as $g \in C$ of 
the cohomology groups of~$\Gamma_g$, and 
the cohomology groups
of~$\Gamma_0(\p \n)$ get replaced similarly by 
taking $R = R_0(\p \n)$
(these replacements are especially needed to have the analog
of Corollary~14(2) of~\cite{greenberg},
where the narrow class number assumption was used implicitly).
%

As usual, let $K_\p$ denote the completion of~$K$ at~$\p$.
With the changes above, the construction of Greenberg goes through mutatis mutandis
to give a point in~$E(K_\p)$ that we again denote~$P_{q,b}$ 
(it is denoted~$P_\psi$ in~\cite{greenberg}, where $b=1$).
Note that while we are using the same notation~$P_{q,b}$ as in Section~\ref{sec:gartner},
it should be clear from the context which point we mean depending on which
construction is used. 
We remark that for his construction,
Greenberg assumes an analog of the conjecture 
of Mazur-Tate-Teitelbaum (conjecture~2 on p.~570 of loc. cit.), and we have
to do the same.
%
%
Just as in~\cite{greenberg}, one has an action 
of~$\Pic(\OO)^+$ on optimal embeddings of~$\OO$ in~$b \widehat{R}^\times b^{-1}$,
and thus on~$P_{q,b}$; we denote this action by~$*$ again.

\begin{conj} \label{conj:greenberg}
The point~$P_{q,b}$ is defined over the narrow ring
class field extension~$K_\OO^+$ of~$K$ associated to the order~$\OO$.
If $\alpha \in \Pic(\OO)^+$, then 
$\rec(\alpha)(P_{q,b}) = \alpha * P_{q,b}$, where 
$\rec: \Pic^+ \OO \ra {\rm Gal}(K_\OO^+/K)$ is the reciprocity
isomorphism of class field theory, as before. 
\end{conj}

Note the similarity of the conjecture above to Conjecture~\ref{conj:gartner}.

\section{Choosing a suitable quaternion algebra}
\label{shdetails}

\comment{
We now discuss more details of the construction 
mentioned in the previous section.
The main point (which we make below) is that 
one can always choose the quaternion algebra $B$ so that one can employ either
the technique
of Gartner or the technique of Greenberg to construct Stark
Heegner points (we will say more
about these two techniques below).
}



Let G1 denote the set of assumptions A, B1, B2, B3, and~B4, and let
G2 denote the set of assumptions A, C1, C2, C3, and C4. If G1
is satisfies, we can carry out 
the construction of~Gartner (as described in 
Section~\ref{sec:gartner}); if G2 holds, 
the construction of~Greenberg (as described in 
Section~\ref{sec:greenberg}) works.

\begin{thm}\label{thm:main}
(i) Suppose that $N$ is square-free. If 
either G1 or G2 hold, then 
the sign in the functional equation of~$E$
over~$K$ 
is~$-1$. \\
(ii) Suppose that the sign in the functional equation 
of~$E$ over~$K$ is~$-1$ and 
the part of~$N$
divisible by primes that are inert in~$K$ is square-free. Then:\\
(a) If there is an archimedian place of~$F$ that is inert in~$K$ (i.e., $K$
is not totally real), then 
one can find a quaternion algebra~$B$ and an Eichler order~$R$ such
G1 holds, i.e., 
one can carry out the construction of~Gartner (as described in 
Section~\ref{sec:gartner}; assuming the conjectures made in the construction).\\
(b) If there is a prime dividing~$N$ that is inert in~$K$, then 
one can find a quaternion algebra~$B$ and an Eichler order~$R$ such
G2 holds, i.e., 
one can carry out the construction 
of~Greenberg (as described in Section~\ref{sec:greenberg}; 
assuming the conjectures made in the construction).\\
(c) One can find a quaternion algebra~$B$ and an Eichler order~$R$ such
that either G1 or G2 hold, i.e.,
one can carry out either the construction of~Gartner (as described in 
Section~\ref{sec:gartner}) 
or the construction of~Greenberg (as described in Section~\ref{sec:greenberg}) 
to construct a Darmon point (assuming the conjectures made in the constructions).
\end{thm}

Note that if
the sign in the functional equation is~$-1$, then
the Birch and Swinnerton-Dyer 
conjecture predicts that $\text{rank } E(K)\ge 1$.  If the
rank is exactly 1, the Gross-Zagier formula would lead one to
expect that (the trace
to $E(K)$ of) the Darmon point has infinite order.  


In the rest of this section, we shall prove Theorem~\ref{thm:main}.
To carry out the construction of Gartner or Greenberg,
we need to find a suitable quaternion algebra~$B$ and an Eichler order~$R$ 
so that all the assumptions made in the construction are satisfied.
We first list the restrictions, and then show when they can be met.
The requirements are as follows:\\

\noi (i) Suppose that Assumption~A holds: there is an embedding of $K$
in~$B$.  This happens if and only if each ramified place of $B$ is 
inert in~$K$.\com{for arch places: if at each place which splits
in~$K$, $B$ is also split.
Thus the ramified places of~$B$ is a 
subset of the set of places of~$F$ that are inert in~$K$. }
Let $\rram$ denote the number of 
real places of~$F$ where $B$ ramifies
and 
$\rinert$ the number
of real places where $B$ is split but that
are inert in~$K$.
The subscript thus indicates which of $B$ and $K$ is non-split, with the understanding that 
we write $\rram$ instead of $r_{\scriptscriptstyle{B,K}}$ since $B$ being non-split implies $K$ is non-split.  \\

\noi (ii) 
Given $K$, the quantity $\rram + \rinert$ is decided,
since it is the number
of real places of~$F$ that are not split
in~$K$. \\

\noi (iii) In Gartner's construction, 
Assumption~B1 says that $\rinert = 1$,
while for Greenberg's construction, 
Assumption~C1 says that $\rinert = 0$ (see the statements just before
the statement of Assumption~C1). \\

In either construction, one needs an Eichler order $R\subset B$ in order to apply the Jacquet-Langlands correspondence.
Let $N^-$ denote the (finite part of the) discriminant of~$B$, 
which is square-free by definition. 
Let $N^+$ and $N'$ be ideals of $\OO_F$ such that $N^-$, $N'$ and
$N^+$ are pairwise relatively prime.  Let $R$ be the Eichler order of $B$ of level $N^+N'$, and put
$\Gamma_0^B(N^+N')={\rm ker}(n: R^\times \ra F^\times)$.
The Jacquet-Langlands correspondence
then says that $S_2(\Gamma_0^{B}(N^+ N'))=S_2(\Gamma_0(N^+ N'N^-))^{N^- -new}$ 
as modules over the Hecke algebra $\mathbb T=\C[\{T_\ell\}_{\ell\nmid N^+ N' N^-}, \{U_p\}_{p|N^+}]$ (the indices~$\ell$ and~$p$ are ideals in~$F$).  
The form associated to~$E$ is in $S_2(\Gamma_0(N))$, so we get the
following conditions on the level of the modular form and the
discriminant of $B$:\\

\noi (iv)  There is a factorization $N=N^+N'N^-$ into three pairwise
coprime ideals.  Here $N^-$ is square-free and divisible only by
primes which are inert in $K$ (as $N^-$ is to serve as the
discriminant ideal of $B$, the second assumption is necessary to
satisfy (i)).  In Greenberg's construction,  $N^-$ is the part of~$N$
divisible by primes other than~$\pp$ that are inert in~$K$, and we
take $N'=\pp$.  Such a factorization exists by Assumption~C3.  In
Gartner's construction, $N^-$ is the product of all prime divisors
of $N$ that are inert in $K$, and $N'=1$. 
\\

\noi (v) By~(i) and~(iv), one sees that the (finite) primes 
where~$B$ is allowed to ramify divide~$N$.
Let $\fram$ denote
the number of primes dividing $N$ 
where $B$ is ramified
and $\finert$ the number of primes dividing $N$ where $B$ splits but that
are inert in~$K$.  
Similar to~(ii), given $K$ and~$N$, the quantity $\fram + \finert$
is independent of $B$, since it is
the number of primes of~$F$ dividing~$N$ that are inert
in~$K$ (recall that we are assuming that $N$ is coprime to the
discriminant of~$K$, so no prime dividing~$N$ ramifies in~$K$). \\
 
\noi (vi) In Gartner's construction, by the extra Assumption~B4,
$\finert = 0$, while in Greenberg's construction,
$\finert = 1$ by Assumptions~C2 and C3.\\

\noi (vii) 
The total number of places where $B$ ramifies is even, so 
$\rram + \fram$ has to be even. And conversely, if 
$\rram + \fram$ is even then a~$B$ exists (ignoring the other conditions). \\

\noi (viii) 
Let $R_b = B \cap b \widehat{R} b^{-1}$.
One needs the existence of an optimal embedding $\OO\hookrightarrow R_b$
%
%
(Assumption~B3 for Gartner's construction and Assumption~C4 in Greenberg's construction).
Such an embedding exists if and only if it exists 
everywhere locally 
(\cite{vig} III.5.11), which happens if and only if 
all the primes 
dividing $N^-$ are inert in $K$ (\cite{vig} II.1.9), and all the primes 
dividing~$N^+$ are split in $K$ (\cite{vig} sentence after II.3.2). 
In (iv), we already had the requirement that 
all the primes dividing $N^-$ are inert, so the only new requirement
is that all the primes 
dividing~$N^+$ are split.
%
%
This requirement is already met in Greenberg's
construction (see~(iv) and Assumption~C3).\\

We now prove part~(i) of Theorem~\ref{thm:main}.
Combining (iii) and~(vi), 
in either construction, $\rinert + \finert =1$, which
combined with~(vii) 
implies that $\rram+ \rinert + \fram + \finert$ is
odd. But $\rram+ \rinert + \fram + \finert$ is precisely 
the total number of places of $K$
where $E$ has a Weierstrass or Tate parametrization,
which in turn is
the exponent of~$-1$ in the sign of the functional equation of 
the $L$-function of~$E$. Thus 
the sign in the functional equation has to be~$-1$.
This proves part~(i) of Theorem~\ref{thm:main}. 

\comment{
The summary of the discussion above is that in order to find a triple
($B$, $R$, $S$) such that 
at least one of Gartner's or Greenberg's construction can be carried out,
one needs the sign in the functional equation to be~$-1$ and 
that the part of~$N$ that is
divisible by primes that are inert in~$K$ is square-free.
}

We next prove part~(ii) of Theorem~\ref{thm:main}.
We shall give two proofs of part~(c). In the first proof, we first try to see 
if the assumptions for Gartner's
construction are satisfied,  and if not, we show that 
the assumptions for Greenberg's
construction hold. In the second proof, we reverse the process:
we first try to satisfy the assumptions for Greenberg's
construction,  and if we can't, we show that 
the assumptions for Gartner's
construction are satisfied. In the process of proving part~(c), we will
prove parts~(a) and~(b).

\begin{proof} (Proof 1 of part~(c) and proof of part~(a))
We first show that 
if $K$ is not totally real, then we can apply Gartner's construction;
this will prove part~(a).
We start with the set of all $B$'s and $R$'s and we
will impose restrictions on this set to satisfy~(i)--(viii)
(for Gartner's construction). The main point 
is that as we impose the restrictions one by one,
at each stage, there should be a choice of $B$ and~$R$ left.
Most of the restrictions in (i)--(viii) are about
ramification of~$B$ at various places. Now a quaternion algebra
with specified ramifications at different places exists if and only if
the number of places where it is ramified is even, which is
condition~(vii).
We impose~(i), and since (vii) can be satisfied while (i) holds, 
we have quaternion algebras~$B$ satisfying~(i)
(this sort of
argument will be used over and over again below, so we will not repeat
the justification we gave in this sentence).
Now the number of
real places of~$F$ that are inert in~$K$ is
$\rram + \rinert$, and so 
$\rram + \rinert$ is non-zero by our hypothesis. While $\rram + \rinert$
is decided by~(ii) (independent of the~$B$'s),
we can restrict to~$B$'s such that $\rinert = 1$ (so that
(iii) is satisfied) and let
$\rram$ be decided by~(ii).
Next we restrict to the $B$'s for which $\finert = 0$
(so that (vi) is satisfied) 
and let $\fram$ be whatever it has to be according to~(v);
however, at this point, we have to check~(vii): we cannot
choose~$\rram$ and~$\fram$ both freely since their sum has to be even.
Now $\rinert + \rram + \finert + \fram = 1 + \rram + 0 + \fram$ is odd
(this parity depends only on the sign of the functional equation of
$L(E_{/K},s)$, hence is independent of the $B$'s), so 
$\rram + \fram$ is even, and (vii) is satisfied, and we are OK.
Thus we can take $N^-$ to be the part of~$N$
divisible by primes that are inert in~$K$ and restrict to those~$B$'s for which
the nonarchimedian places where $B$ ramifies are precisely the ones dividing~$N^-$
(here we are using the hypothesis that 
the part of~$N$ divisible by primes that are inert in~$K$ is square-free).
We take $N^+ = N/N^-$, so that (iv) and (viii) are satisfied
(note that $N' = 1$) and choose an order~$R$ of level~$N^+$.
Thus we can find a~$B$ and an~$R$
for which (i)--(viii) are satisfied for Gartner's construction.

If $K$ is totally real,
then we claim that we can apply 
Greenberg's construction.
We start with the set of all $B$'s and $R$'s and impose~(i).
By (ii), $\rinert + \rram = 0$, so $\rinert = \rram = 0$, and (iii) is satisfied
(for Greenberg's construction).
Now $\finert + \fram = \rinert + \rram + \finert + \fram$ is odd, 
and in particular, non-zero. So 
we may restrict to~$B$'s such that $\finert = 1$ (then
(vi) is satisfied) and let $\fram$ be whatever it needs to be
to satisfy~(v) ($\fram$ will be even).
Again, at this point, we have
to check that there are $B$'s left satisfying the conditions above
since by~(vii), $\rram + \fram$ has to be even; but this is true
since $\rram=0$ and~$\fram$ is even as mentioned above.
Thus we take a prime $\pp$ that divides~$N$ and is inert
in~$K$, and let $N^-$ be the product of all
primes except~$\pp$ that divide~$N$ and are inert
in~$K$ (here we are using the hypothesis that the part of~$N$ 
divisible by primes that are inert in~$K$ is square-free).
We take $N' = \pp$, $N^+ = N/ (N^- \pp)$ so that (viii) is satisfied.
Also, as mentioned above, (iv) is automatic for Greenberg's construction.
Thus we can find a~$B$ and an~$R$
for which (i)--(viii) are satisfied for Greenberg's construction.
\end{proof}

\begin{proof} (Proof 2 of part~(c)  and proof of part~(b))
We first show that 
if $N$ is divisible by a prime that is
inert in~$K$, then we can apply Greenberg's construction;
this will prove part~(b).
As in Proof~1, we start with the set of all $B$'s and $R$'s and we
impose restrictions on this set to satisfy~(i)--(viii) 
(for Greenberg's construction). The main point 
is that as we impose the restrictions one by one,
at each stage, there should be a choice of $B$ and~$R$ left.
Most of the restrictions in (i)--(viii) are about
ramification of~$B$ at various places. Now a quaternion algebra
with specified ramifications exists if and only if
the number of places where it is ramified is even, which is
condition~(vii).
We impose~(i), and since (vii) can be satisfied while (i) holds, 
we have quaternion algebras~$B$ satisfying~(i)
We pick a prime~$\p$ such that $\p | N$ and $\p$ is inert in~$K$.
We restrict to the
$B$'s such that the (finite) primes where $B$ is ramified
is precisely  the set of primes  except~$\p$ that 
divide~$N$ and are inert in~$K$. Thus
$N^-$ is the product of all
primes except for~$\p$ that divide~$N$ and are inert in~$K$.
(so (v) is satisfied).
We take 
$N^+ = N/ N^- \p$ and $N' = \pp$; then (iv), (vi), and~(viii) are satisfied
(here we are using the hypothesis that 
the part of~$N$ that is
divisible by primes that are inert in~$K$ is square-free).
We further restrict to the $B$'s such that $\rinert = 0$ (so that (iii)
is satisfied) and let $\rram$ be decided by~(ii);
however, at this point, we have to check~(vii): we cannot
choose~$\rram$ and~$\fram$ both freely since their sum has to be even.
Now $\rinert + \rram + \finert + \fram = 0 + \rram + 1 + \fram$ is odd
(this parity depends only on~$E$ and~$K$, and is independent of the $B$'s), so 
$\rram + \fram$ is even, and (vii) is satisfied.
Thus (i)--(viii) are satisfied for Greenberg's construction.

If the part of~$N$
divisible by primes that are inert in~$K$ is empty, then by~(v),
$\finert = \fram = 0$ (so (vi) is satisfied for
Gartner's construction), and
we restrict
to~$B$'s that are not ramified at any (finite) prime (so $N^- = 1$).
Now $\rinert + \rram = \rinert + \rram + \finert + \fram $ is odd
by hypothesis, hence non-zero.
We restrict to $B$'s such that $\rinert = 1$ (so that (iii) is
satisfied) and let $\rram$ be decided by~(ii)  ($\rram$ will be even). 
Again, at this point, we have
to check that there are $B$'s left satisfying the conditions above
since by~(vii), $\rram + \fram$ has to be even; but this is true
since $\fram=0$ and~$\rram$ is even as mentioned above.
We then take $N^- = N' = 1$, and $N^+ = N$, so that
(iv) and (viii) are satified.  
Thus (i)--(viii) are satisfied for Gartner's construction.
\end{proof}

\begin{rmk} \label{rmk:assn}
(i) We made Asssumption~B4 in Gartner's construction (Section~\ref{sec:gartner})
in order to get part~(i) of Theorem~\ref{thm:main}.
Also, this assumption is a natural choice to be made in the construction anyway.
If Assumption~B4 is dropped, then 
part~(ii) of Theorem~\ref{thm:main}
is still true, and in particular, if
the sign in the functional equation is~$-1$ and 
the part of~$N$
divisible by primes that are inert in~$K$ is square-free, then
one can carry out either the construction of~Gartner 
or the construction of~Greenberg to construct a Darmon point.\\
(ii) There were choices for the quaternion algebra~$B$ and
the Eichler order~$R$
in what we did above for the proof of part~(ii) 
and there may be
other ways of applying Greenberg's or Gartner's constructions than what
we did. Also,  
if the sign in the functional equation is~$-1$ and 
the part of~$N$
divisible by primes that are inert in~$K$ is square-free,
and there is at least one real place and one prime dividing~$N$ that
are inert in~$K$, then either of
Greenberg's or Gartner's constructions can be carried out
(by the first paragraphs of Proofs~1 and~2).
It would be interesting to see if and how the Darmon points one
gets by different choices (when available) are related.
\end{rmk}


\providecommand{\bysame}{\leavevmode\hbox to3em{\hrulefill}\thinspace}
\providecommand{\MR}{\relax\ifhmode\unskip\space\fi MR }
\providecommand{\MRhref}[2]{%
  \href{http://www.ams.org/mathscinet-getitem?mr=#1}{#2}
}
\providecommand{\href}[2]{#2}

\end{document}